\documentclass[12pt]{article}

\usepackage[margin=1in]{geometry}
\usepackage{amsmath,amssymb}
\usepackage{graphicx}
\usepackage[colorlinks=true,linkcolor=blue,citecolor=blue,urlcolor=blue]{hyperref}
\usepackage{natbib}
\usepackage{booktabs}
\usepackage{caption}
\usepackage{subcaption}
\usepackage{xcolor}
\usepackage{setspace}
\usepackage{enumitem}

\title{\textbf{Towards Sustainable Hydrogen Systems: Supply Chain Optimization
  with Model Predictive Control and Reinforcement Learning}}

\author{Mahammad Valiyev\\[4pt]
  \small Department of Chemical Engineering\\
  \small University of Southern California\\
  \small Los Angeles, California, United States}

\date{}

\begin{document}
\maketitle

\begin{abstract}
Hydrogen supply chains are expected to play a central role in future low-carbon
energy systems. However, the operation of such systems is inherently challenged
by factors such as uncertainty in renewable availability, electricity prices, and
hydrogen demand, as well as by practical engineering constraints associated with
electrolyzers, energy storage, and grid interaction. Effective control strategies
must therefore balance economic performance, operational feasibility, and
sustainability objectives under time-varying and uncertain conditions, especially
as hydrogen infrastructure scales beyond pilot deployments.

This paper investigates and compares four operational control approaches for a
renewable-powered hydrogen supply chain: a rule-based controller (RBC), model
predictive control (MPC), reinforcement learning without forecasts (RL-NF), and
reinforcement learning with forecast-augmented observations (RL-F). All methods
are evaluated within a unified, physically realistic framework that incorporates
electrolyzer minimum load and ramping constraints, battery and hydrogen storage
dynamics, grid import limits, renewable curtailment, and consistent economic
assumptions. This unified formulation ensures methodological consistency, enabling
a fair and rigorous comparison of predictive and adaptive control paradigms under
identical system conditions and uncertainty realizations.

Simulation results show that MPC achieves the highest economic performance by
effectively exploiting short-term forecasts to reduce grid dependence and optimize
storage utilization, highlighting the value of predictive optimization when
reliable forecasts are available. RL-NF demonstrates robust and competitive
performance despite operating without future information, while RL-F does not
consistently outperform its no-forecast counterpart, illustrating the challenges
of learning from noisy forecasts in complex energy systems. These findings
illustrate the relative strengths and limitations of predictive and learning-based
control approaches for hydrogen supply chains under uncertainty and motivate
further investigation into their roles in large-scale hydrogen system operation.
\end{abstract}

\section{Introduction}
\label{sec:intro}

Global decarbonization agendas increasingly emphasize the need for scalable,
low-carbon energy carriers capable of integrating high shares of renewable
generation. In this context, hydrogen has emerged as a promising option for
long-duration energy storage, industrial decarbonization, and clean fuel supply
for transportation and power systems \citep{IEA2019,Staffell2019}. When produced
using renewable electricity, hydrogen can enable deep emissions reductions while
enhancing system flexibility in energy systems with growing shares of variable
wind and solar generation. However, realizing this potential at scale requires not
only infrastructure investment, but also effective operational strategies capable
of managing uncertainty, variability, and complex physical constraints.

The operation of hydrogen supply chains presents significant challenges due to
uncertainty in renewable availability, electricity market prices, and hydrogen
demand, as well as engineering constraints associated with electrolyzers, energy
storage systems, and grid interaction. Electrolyzers exhibit minimum loading and
ramping constraints, energy storage assets are subject to capacity and throughput
limits, and grid imports are often constrained by economic and regulatory
considerations. As hydrogen systems scale beyond pilot projects, suboptimal
operational decisions can lead to excessive grid dependence, renewable curtailment,
increased costs, or underutilized assets. Consequently, optimizing hydrogen supply
chain operation under uncertainty has become a critical research problem at the
intersection of energy systems engineering, optimization, and control.

A wide range of control and optimization approaches have been explored for energy
systems and energy supply chains. Rule-based control strategies remain common in
practice due to their simplicity, interpretability, and ease of implementation,
particularly in early-stage or pilot-scale deployments
\citep{Oldewurtel2012,Behzadi2023}. However, such approaches are inherently
reactive and typically struggle to exploit future information or adapt efficiently
to changing operating conditions as system complexity increases. Model Predictive
Control (MPC) has therefore received significant attention in recent years for
hydrogen-integrated energy systems, including renewable microgrids,
power-to-hydrogen facilities, and hybrid storage systems, owing to its ability to
explicitly incorporate forecasts, operational constraints, and economic objectives
within a receding-horizon optimization framework
\citep{Parisio2014,Rawlings2020,Brka2016,Velarde2017}. Recent studies demonstrate
that MPC can substantially reduce grid reliance, renewable curtailment, and
operating costs in hydrogen-coupled systems when short-term forecasts are
sufficiently reliable, though at the expense of increased computational burden and
sensitivity to forecast errors \citep{Huang2022}.

More recently, reinforcement learning (RL) has emerged as a data-driven
alternative for controlling complex energy and hydrogen systems under uncertainty,
particularly in settings where accurate system models or reliable forecasts are
difficult to obtain \citep{Glavic2017,VazquezCanteli2019}. Advances in deep RL
and safe learning frameworks have enabled applications to battery dispatch,
integrated energy systems, and hydrogen production scheduling
\citep{Ruelens2016,Yang2021,Liang2024,Ye2024}. These studies suggest that RL can
learn adaptive policies capable of handling nonlinear dynamics and stochastic
disturbances, but also highlight persistent challenges related to constraint
satisfaction, training stability, and interpretability. In addition, recent work
indicates that augmenting RL policies with forecast information does not
necessarily guarantee improved performance, particularly when forecasts are noisy
or only weakly informative \citep{Oh2020}. As a result, it remains an open
question under which conditions forecast-augmented RL approaches can outperform
reactive RL policies or established predictive optimization methods such as MPC
in hydrogen supply-chain operations.

This paper contributes to the ongoing discussion by presenting a unified
comparative study of four operational control strategies for a renewable-powered
hydrogen supply chain: RBC, MPC, RL-NF, and RL-F. The modeled system captures
the core elements of emerging green hydrogen configurations, including variable
renewable generation, battery storage, an electrolyzer with minimum load and
ramping constraints, hydrogen storage with throughput limits, grid interaction,
and uncertain demand. This unified formulation enables a rigorous apples-to-apples
comparison of predictive and learning-based control paradigms under identical
disturbances, constraints, and economic assumptions.

The results demonstrate clear differences in how predictive and adaptive control
paradigms manage uncertainty and system constraints. MPC achieves the highest
economic performance by effectively exploiting short-term forecasts to reduce
grid reliance and optimize storage utilization. RL-NF exhibits robust and
competitive performance despite operating under partial observability, learning
conservative policies that prioritize renewable utilization and high hydrogen
utilization rates. In contrast, RL-F does not consistently outperform its
no-forecast counterpart, highlighting the challenges associated with learning
from noisy forecast information in complex, constrained energy systems. These
findings motivate further investigation into robust and hybrid control strategies
as hydrogen infrastructure scales.

\section{Methodology}
\label{sec:methodology}

\subsection{System Description and Problem Formulation}
\label{subsec:system}

\subsubsection{Hydrogen Supply Chain Process Model}

\begin{figure}[htbp]
  \centering
  \includegraphics[width=0.92\textwidth]{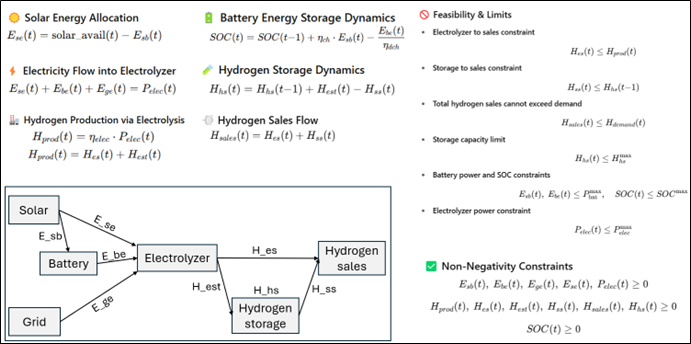}
  \caption{Process model, mass and energy balance, feasibility constraints, and
    operational limits for the renewable-powered hydrogen supply chain. The diagram
    links solar generation, battery storage, grid interaction, electrolyzer
    operation, hydrogen storage, and sales through energy balance relationships that
    define the feasible control space.}
  \label{fig:process_model}
\end{figure}

This study considers a stylized renewable-powered hydrogen supply chain designed
to capture the key operational interactions between variable renewable generation,
energy storage, hydrogen production, and market demand. The system consists of
four main components illustrated in Figure~\ref{fig:process_model}.

\paragraph{Renewable electricity source.}
Solar photovoltaic generation provides a variable, low-cost electricity input to
the system. Solar availability is time-dependent and exhibits diurnal patterns
with stochastic variability, representing typical photovoltaic generation behavior.

\paragraph{Battery energy storage system (BESS).}
A battery is used to temporarily store electricity by charging during periods of
excess solar generation and discharging when renewable availability is
insufficient. The battery is characterized by finite power and energy capacity
limits, as well as charging and discharging efficiencies.

\paragraph{Electrolyzer.}
The electrolyzer converts electrical energy into hydrogen with a fixed conversion
efficiency. Electricity supplied to the electrolyzer can originate from solar
generation, battery discharge, or the electricity grid. The electrolyzer is
subject to rated power limits, minimum loading constraints, and ramping
constraints.

\paragraph{Hydrogen storage and demand.}
Produced hydrogen can either be sold immediately to meet exogenous demand or
stored in a hydrogen storage tank with limited capacity. Stored hydrogen may later
be withdrawn to satisfy demand shortfalls, subject to storage throughput
constraints.

The system operates in discrete hourly time steps over a one-week simulation
horizon and is exposed to time-varying solar availability, electricity prices, and
hydrogen demand.

\subsection{State Variables, Control Variables, and System Dynamics}
\label{subsec:dynamics}

\subsubsection{State Variables}

System dynamics are governed by two primary state variables:
\begin{itemize}
  \item \textbf{Battery state of charge (SOC):}
    $\mathrm{SOC}_t \in [0,\mathrm{SOC}_{\max}]$,
    representing the stored electrical energy in the battery at time $t$.
  \item \textbf{Hydrogen storage level:}
    $H_t \in [0, H_{\max}]$,
    representing the amount of hydrogen stored in the tank.
\end{itemize}
These states evolve according to energy and mass balance equations driven by
control actions and system efficiencies.

\subsubsection{Control Variables}

At each time step the controller determines the following energy and hydrogen
allocation decisions:
\begin{itemize}[leftmargin=2em]
  \item $E_{sb,t}$: electricity charged from solar into the battery
  \item $E_{be,t}$: electricity discharged from the battery to the electrolyzer
  \item $E_{ge,t}$: electricity drawn from the grid to the electrolyzer
  \item $E_{se,t}$: electricity supplied directly from solar to the electrolyzer
  \item $P_{e,t}$: total electrical power supplied to the electrolyzer
  \item $H_{es,t}$: hydrogen sold directly from production
  \item $H_{est,t}$: hydrogen sent to storage
  \item $H_{ss,t}$: hydrogen sold from storage
\end{itemize}
In the MPC formulation these variables are optimized explicitly. In the RL
formulations, high-level agent actions are mapped to feasible values through a
feasibility layer, ensuring consistent system behavior and constraint satisfaction
across all methods.

\subsubsection{System Dynamics}

The battery state of charge evolves according to:
\begin{equation}
  \mathrm{SOC}_{t+1} = \mathrm{SOC}_t + \eta_{ch}\,E_{sb,t}
    - \frac{E_{be,t}}{\eta_{dch}}.
  \label{eq:soc}
\end{equation}

Hydrogen storage dynamics are given by:
\begin{equation}
  H_{t+1} = H_t + H_{est,t} - H_{ss,t}.
  \label{eq:hstorage_dyn}
\end{equation}

Hydrogen production from the electrolyzer is modeled as:
\begin{equation}
  H_{prod,t} = \eta_{elec}\,P_{e,t}.
  \label{eq:hprod}
\end{equation}

Total hydrogen sales satisfy:
\begin{equation}
  H_{sales,t} = H_{es,t} + H_{ss,t}.
  \label{eq:hsales}
\end{equation}

\subsection{Objective Function}
\label{subsec:objective}

The operational objective is to maximize net profit over the simulation horizon,
equivalently minimizing total operating cost minus hydrogen sales revenue. At each
time step the instantaneous cost includes: grid electricity cost (including a
carbon adder), battery cycling cost, electrolyzer operating cost, and hydrogen
storage holding cost. Revenue is generated through hydrogen sales. The per-step
cost is:
\begin{equation}
  \min_{u_t}\;
  \Bigl[(\pi_t + c_{CO_2})\,E_{ge,t}
    + c_{bat}(E_{sb,t} + E_{be,t})
    + c_{elec}\,P_{e,t}
    + c_{store}\,H_t
    - p_H\,H_{sales,t}\Bigr].
  \label{eq:objective}
\end{equation}
This objective is used directly in the MPC formulation and serves as the reward
signal for the RL controllers, enabling consistent economic comparison across
methods.

\subsection{Operational Constraints}
\label{subsec:constraints}

System operation is subject to the following constraints:
\begin{align}
  &0 \leq P_{e,t} \leq P_{\max},\quad
   P_{e,t} \geq P_{\min}\cdot z_t
   \label{eq:elec_power}\\[4pt]
  &0 \leq E_{sb,t},\;E_{be,t} \leq P_{bat,\max},\quad
   0 \leq \mathrm{SOC}_t \leq \mathrm{SOC}_{\max}
   \label{eq:bat}\\[4pt]
  &0 \leq H_t \leq H_{\max}
   \label{eq:hstorage_bounds}\\[4pt]
  &H_{prod,t} = H_{es,t} + H_{est,t}
   \label{eq:hbalance}\\[4pt]
  &P_{e,t} = E_{se,t} + E_{be,t} + E_{ge,t}
   \label{eq:ebalance}\\[4pt]
  &H_{sales,t} \leq D_t
   \label{eq:demand}\\[4pt]
  &E_{se,t} + E_{sb,t} \leq S_t
   \label{eq:solar_limit}
\end{align}
These constraints are enforced explicitly in the MPC formulation and implicitly
through a feasibility layer in the RL controllers, ensuring physically realistic
and safe operation under all control strategies. Unmet hydrogen demand is
permitted and does not incur an explicit penalty beyond foregone sales revenue.

\subsection{Disturbance and Forecast Models}
\label{subsec:disturbances}

The hydrogen supply chain is driven by three exogenous, time-varying disturbances
evaluated at hourly resolution $t = 0,1,\ldots,T-1$. All controllers are
evaluated under identical disturbance realizations to ensure a fair comparison.

\subsubsection{Solar Availability}

\begin{figure}[htbp]
  \centering
  \includegraphics[width=0.92\textwidth]{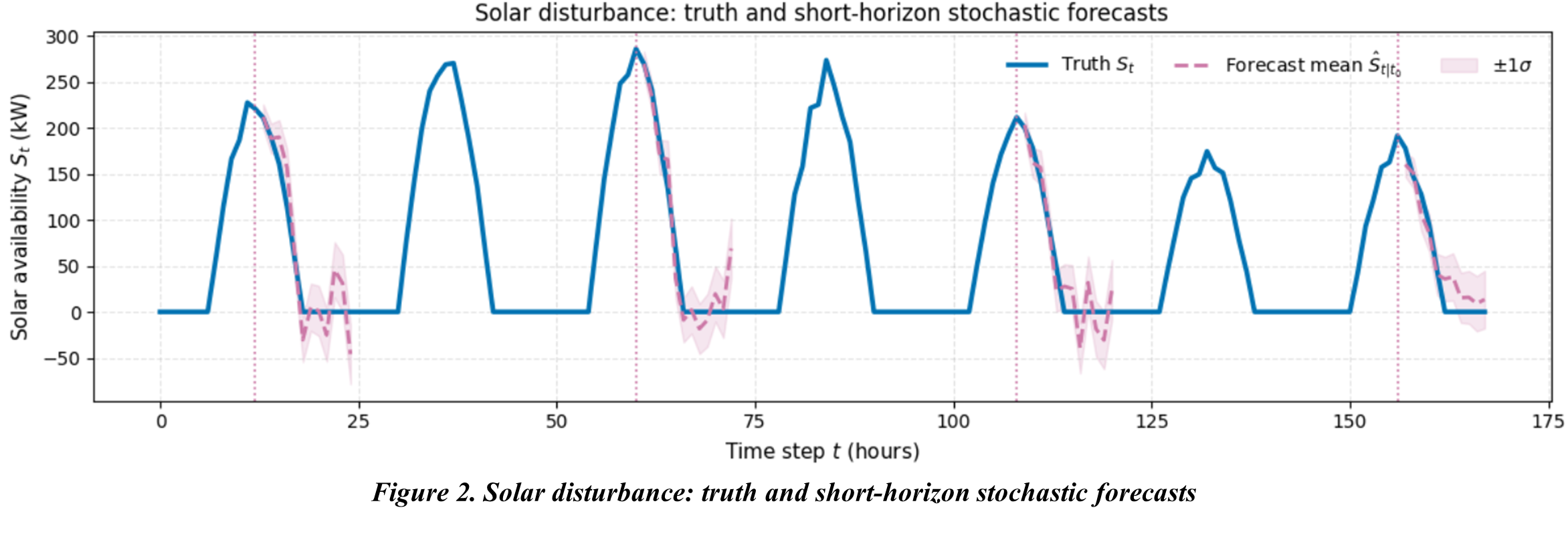}
  \caption{Solar availability disturbance (realized) and representative
    short-horizon stochastic forecasts used by the predictive controllers.}
  \label{fig:solar_forecast}
\end{figure}

Solar power availability $S_t\;[\mathrm{kW}]$ is modeled as:
\begin{equation}
  S_t = \max\!\bigl(0,\;P_{PV}^{\max}\,b(h_t)\,(c(d_t)+\varepsilon_t^{S})\bigr),
  \label{eq:solar}
\end{equation}
where $P_{PV}^{\max}=300\;\mathrm{kW}$ is the installed PV capacity,
$h_t = t \bmod 24$ is the hour of day, and $d_t = \lfloor t/24 \rfloor$ is the
day index. The deterministic diurnal shape is a half-sine function:
\begin{equation}
  b(h) = \begin{cases}
    0, & h < 6 \;\text{or}\; h > 18,\\[2pt]
    \displaystyle\sin\!\left(\pi\,\frac{h-6}{12}\right), & 6 \leq h \leq 18,
  \end{cases}
  \label{eq:solar_shape}
\end{equation}
while slow cloud-induced variability is modeled as:
\begin{equation}
  c(d) = 0.75 + 0.20\sin\!\left(\frac{2\pi d}{7}\right).
  \label{eq:cloud}
\end{equation}
Small, zero-mean Gaussian noise $\varepsilon_t^{S}\sim\mathcal{N}(0,0.03^2)$
captures short-term irradiance fluctuations. The realized solar disturbance and
representative short-horizon stochastic forecasts are illustrated in
Figure~\ref{fig:solar_forecast}.

\subsubsection{Electricity Prices}

\begin{figure}[htbp]
  \centering
  \includegraphics[width=0.92\textwidth]{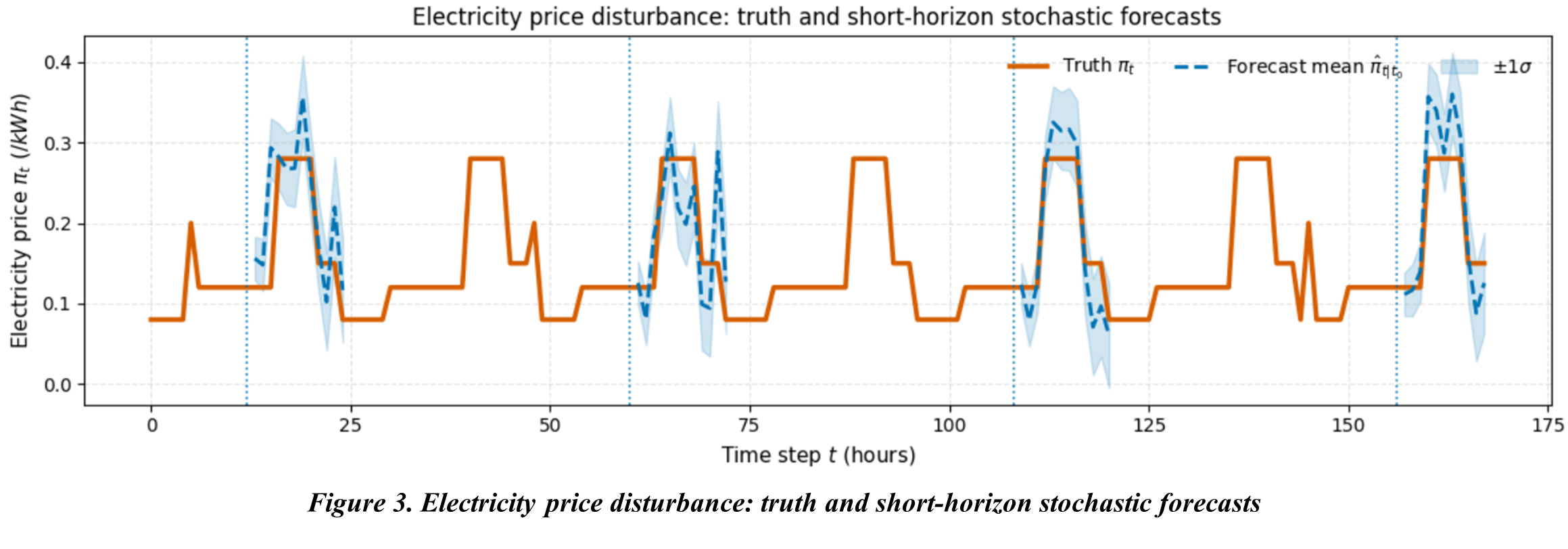}
  \caption{Electricity price disturbance (realized) and representative
    short-horizon stochastic forecasts.}
  \label{fig:price_forecast}
\end{figure}

The grid electricity price $\pi_t\;[\$/\mathrm{kWh}]$ follows a time-of-use
(TOU) tariff structure with occasional price spikes:
\begin{equation}
  \pi_t = \pi_{TOU}(h_t)\cdot\xi_t,
  \label{eq:price}
\end{equation}
where the baseline TOU price is:
\begin{equation}
  \pi_{TOU}(h) = \begin{cases}
    0.08, & 0 \leq h < 6,\\
    0.12, & 6 \leq h < 16,\\
    0.28, & 16 \leq h < 21,\\
    0.15, & 21 \leq h < 24,
  \end{cases}
  \label{eq:tou}
\end{equation}
and the stochastic spike multiplier is:
\begin{equation}
  \xi_t = \begin{cases}
    2.5, & \text{with probability } 0.02,\\
    1,   & \text{otherwise.}
  \end{cases}
  \label{eq:spike}
\end{equation}
The realized electricity price disturbance and representative forecasts are shown
in Figure~\ref{fig:price_forecast}.

\subsubsection{Hydrogen Demand}

\begin{figure}[htbp]
  \centering
  \includegraphics[width=0.92\textwidth]{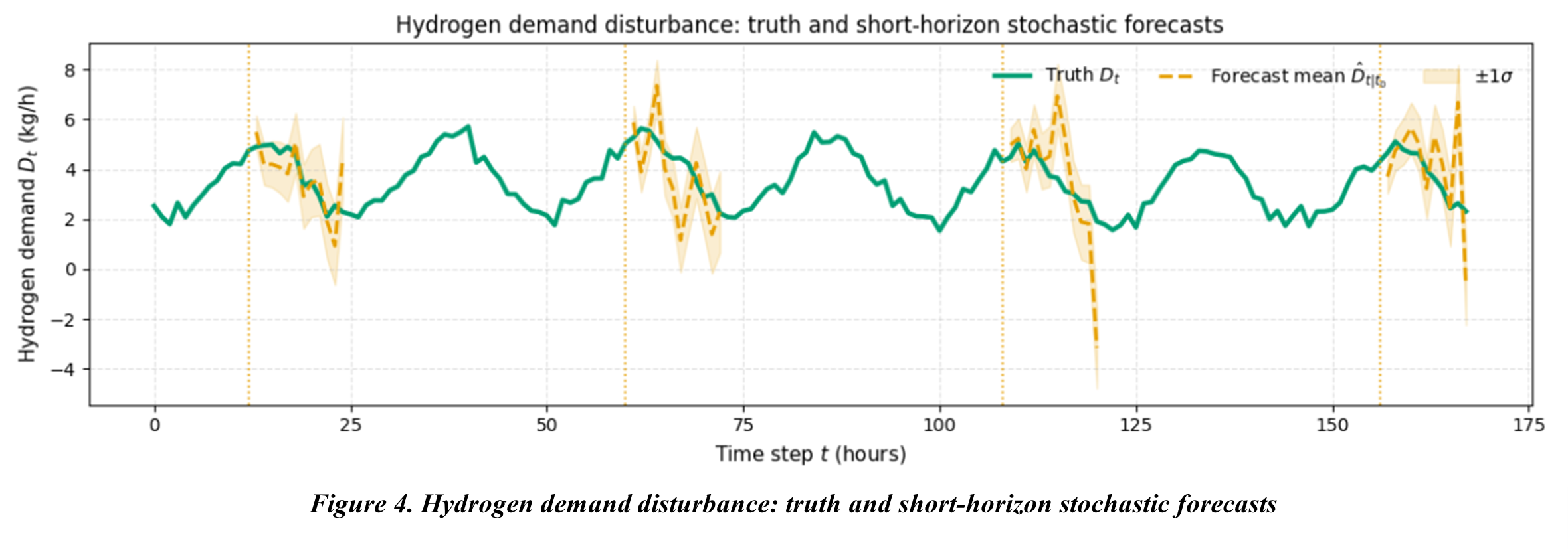}
  \caption{Hydrogen demand disturbance (realized) and representative short-horizon
    stochastic forecasts.}
  \label{fig:demand_forecast}
\end{figure}

Hydrogen demand $D_t\;[\mathrm{kg/h}]$ is modeled as:
\begin{equation}
  D_t = \max\!\left(0,\;3.5
    + 1.5\sin\!\left(\frac{2\pi(h_t-8)}{24}\right)
    + 0.3\sin\!\left(\frac{2\pi d_t}{7}\right)
    + \varepsilon_t^{D}\right),
  \label{eq:demand_model}
\end{equation}
where $\varepsilon_t^{D}\sim\mathcal{N}(0,0.25^2)$. The demand magnitude is
intentionally scaled to be consistent with a 300~kW electrolyzer, corresponding
to a maximum production rate of approximately 6~kg/h given a conversion factor of
50~kWh/kg. The realized hydrogen demand disturbance and representative forecasts
are shown in Figure~\ref{fig:demand_forecast}.

\subsubsection{Forecast Availability and Uncertainty}

Disturbances are assumed to be known perfectly only at the current time step.
Short-horizon stochastic forecasts are available to the MPC controller and RL-F,
while RL-NF operates under partial observability. For any disturbance signal
$x_t\in\{S_t,\pi_t,D_t\}$, a forecast issued at time $t_0$ for lead time
$\delta=1,\ldots,H$ is modeled as:
\begin{equation}
  \hat{x}_{t_0+\delta\mid t_0} = x_{t_0+\delta} + \eta_\delta,\quad
  \eta_\delta\sim\mathcal{N}\!\left(0,\sigma_x^2(\delta)\right),
  \label{eq:forecast_model}
\end{equation}
where $H=12$ hours is the prediction horizon. Forecast uncertainty grows with
lead time according to:
\begin{equation}
  \sigma_x(\delta) = \sigma_{x,0}\sqrt{1+\alpha\delta}.
  \label{eq:forecast_uncertainty}
\end{equation}
The base uncertainty levels are signal-specific:
$\sigma_{S,0}=10\;\mathrm{kW}$,
$\sigma_{\pi,0}=0.02\;\$/\mathrm{kWh}$,
$\sigma_{D,0}=0.5\;\mathrm{kg/h}$,
and $\alpha=0.8$ controls the rate of uncertainty growth. All forecasts are
generated using identical disturbance models and random seeds across control
strategies, ensuring that performance differences arise solely from control logic
rather than informational asymmetry.

\subsection{Control Strategies}
\label{subsec:control}

\begin{figure}[htbp]
  \centering
  \includegraphics[width=0.92\textwidth]{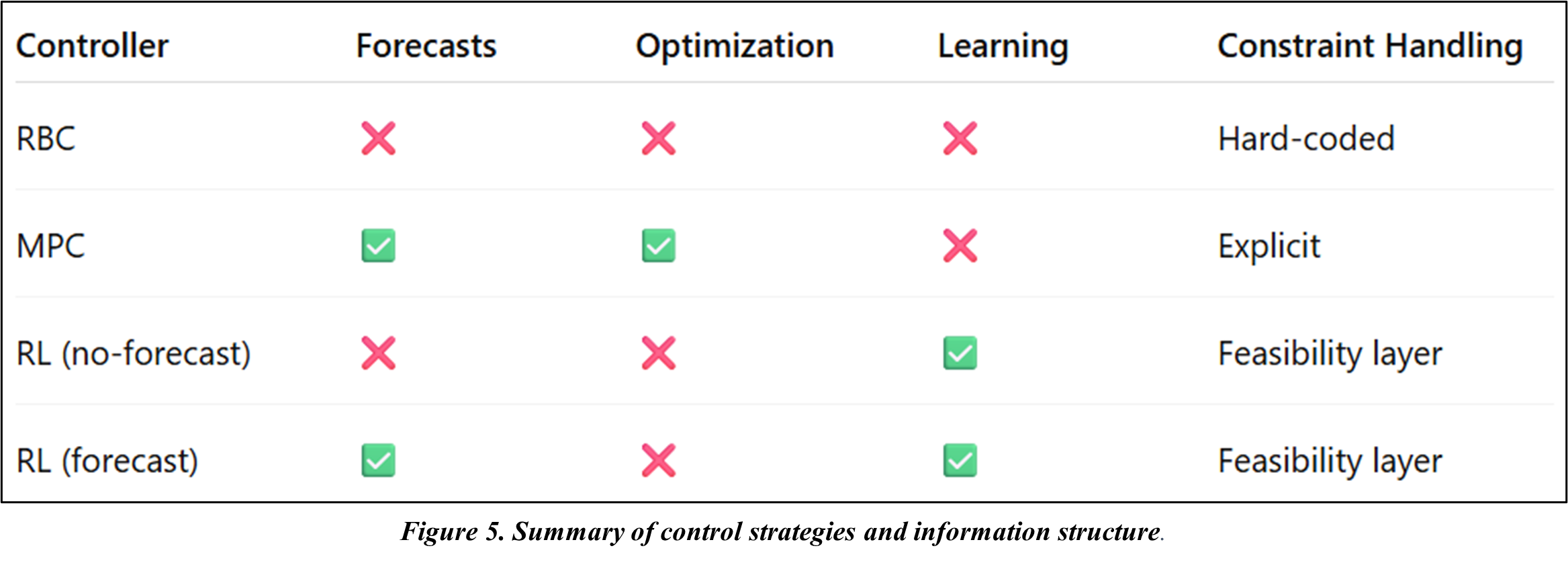}
  \caption{Summary of the four control strategies and their information
    structures. RBC relies on hard-coded rules without forecasts or learning; MPC
    leverages forecasts with explicit optimization; RL-NF and RL-F rely on
    learning, with RL-F incorporating short-horizon predictions. Constraint
    enforcement differs: MPC uses explicit constraints; RL methods employ a
    feasibility layer.}
  \label{fig:control_summary}
\end{figure}

This section presents the four control strategies evaluated in this study. All
controllers operate on the same system model, constraints, cost structure, and
disturbance realizations. Differences arise solely from how control decisions are
generated — in particular in the use of forecasts, optimization, and learning.
The structural differences are summarized in Figure~\ref{fig:control_summary}.

\subsubsection{Rule-Based Control (RBC)}

Rule-based control (RBC) represents a heuristic, non-optimizing baseline
commonly used in practice due to its simplicity and transparency. The controller
operates myopically, relying only on current system states and disturbances,
without forecasts or explicit optimization.

At each time step $t$ the RBC policy follows a fixed priority structure:
\begin{enumerate}[leftmargin=2em]
  \item Available solar energy is allocated to the electrolyzer first, up to its
    rated capacity $P_{\max}$. Excess solar is used to charge the battery if it is
    not near full.
  \item If solar is insufficient to reach $P_{\max}$ and
    $\mathrm{SOC}_t > 0.2\,\mathrm{SOC}_{\max}$, the battery is discharged to
    supply the electrolyzer.
  \item Any remaining electrolyzer power requirement is met by drawing electricity
    from the grid, independent of electricity price.
  \item Hydrogen produced at time $t$ is first used to meet current demand $D_t$.
    Surplus hydrogen is stored (subject to capacity); if production is insufficient,
    hydrogen is withdrawn from storage.
\end{enumerate}

RBC serves as a benchmark baseline, providing a lower bound on performance against
which more advanced control strategies are evaluated.

\subsubsection{Model Predictive Control (MPC)}

MPC formulates system operation as a rolling-horizon optimization problem that
explicitly accounts for system dynamics, constraints, and forecasts of future
disturbances. At each time step $t$, MPC solves a finite-horizon optimization
problem over prediction horizon $H$, using forecasts $\hat{S}_{t+k}$,
$\hat{\pi}_{t+k}$, and $\hat{D}_{t+k}$ for $k=0,\ldots,H-1$. The optimization
variables are
$\{E_{sb,t+k}, E_{be,t+k}, E_{ge,t+k}, P_{e,t+k},
   H_{es,t+k}, H_{est,t+k}, H_{ss,t+k}\}_{k=0}^{H-1}$,
and state dynamics follow equations \eqref{eq:soc}--\eqref{eq:hprod}.

The MPC minimizes cumulative economic cost over the horizon:
\begin{equation}
  \min\sum_{k=0}^{H-1}\!\Bigl[
    (\hat{\pi}_{t+k}+c_{CO_2})\,E_{ge,t+k}
    + c_{bat}(E_{sb,t+k}+E_{be,t+k})
    + c_{elec}\,P_{e,t+k}
    + c_{store}\,H_{t+k}
    - p_H\,H_{sales,t+k}
  \Bigr].
  \label{eq:mpc_obj}
\end{equation}
All physical and operational constraints \eqref{eq:elec_power}--\eqref{eq:solar_limit}
are enforced explicitly. Only the first control action is implemented; the horizon
then shifts forward by one step, forecasts are updated, and the problem is
re-solved.

\subsubsection{Reinforcement Learning Framework}

RL controllers learn a control policy through interaction with the environment,
without solving an explicit optimization problem at runtime. The hydrogen supply
chain is formulated as a Markov decision process (MDP):
\begin{itemize}
  \item \textbf{State}: system states and exogenous information available to the
    agent.
  \item \textbf{Action}: high-level continuous control vector
    $a_t = [u_{pe,t},\,u_{bat,t},\,u_{store,t}]^{\top}\in[0,1]^3$,
    controlling electrolyzer utilization, battery
    charging/discharging preference, and hydrogen storage preference.
  \item \textbf{Transition}: governed by the physical dynamics
    \eqref{eq:soc}--\eqref{eq:hsales}.
  \item \textbf{Reward}: per-step economic profit
    \begin{equation}
      r_t = p_H H_{sales,t}
        - \bigl[(\pi_t+c_{CO_2})\,E_{ge,t}
        + c_{bat}(E_{sb,t}+E_{be,t})
        + c_{elec}\,P_{e,t}
        + c_{store}\,H_t\bigr],
      \label{eq:rl_reward}
    \end{equation}
    identical to the MPC objective on a per-step basis.
\end{itemize}
All RL agents are trained using Proximal Policy Optimization (PPO) with
continuous action spaces, identical hyperparameters, and identical disturbance
realizations to ensure fair comparison.

\paragraph{RL without forecasts (RL-NF).}
The agent observes only current information:
$[\mathrm{SOC}_t,\,H_t,\,S_t,\,\pi_t,\,D_t]$, normalized to fixed ranges.
Because future disturbances are not included, the agent operates under partial
observability; anticipating future price spikes, solar availability, or demand
changes is not possible, and learned policies tend to be reactive.

\paragraph{RL with forecasts (RL-F).}
The forecast-augmented variant extends the observation to include short-horizon
noisy forecasts:
$[\hat{S}_t,\hat{\pi}_t,\hat{D}_t,\ldots,
  \hat{S}_{t+H_f-1},\hat{\pi}_{t+H_f-1},\hat{D}_{t+H_f-1}]$,
generated using the same stochastic forecast model employed by MPC. Despite
access to forecasts, RL-F differs fundamentally from MPC: it does not solve a
multi-step optimization problem, constraint coupling across time is implicit
rather than explicit, and planning emerges indirectly through policy learning.

\paragraph{Feasibility layer.}
In both RL formulations, actions are interpreted as normalized preferences and
mapped deterministically to admissible control variables through a feasibility
layer that (i) clips control targets to allowable ranges, (ii) allocates
available solar, battery, and grid power sequentially to satisfy electrolyzer
demand, and (iii) ensures energy and mass balance at each time step. This
guarantees physically realistic operation and enables a fair comparison with MPC.

\section{Results}
\label{sec:results}

\subsection{Rule-Based Control: Baseline Performance}
\label{subsec:results_rbc}

\begin{figure}[htbp]
  \centering
  \includegraphics[width=0.92\textwidth]{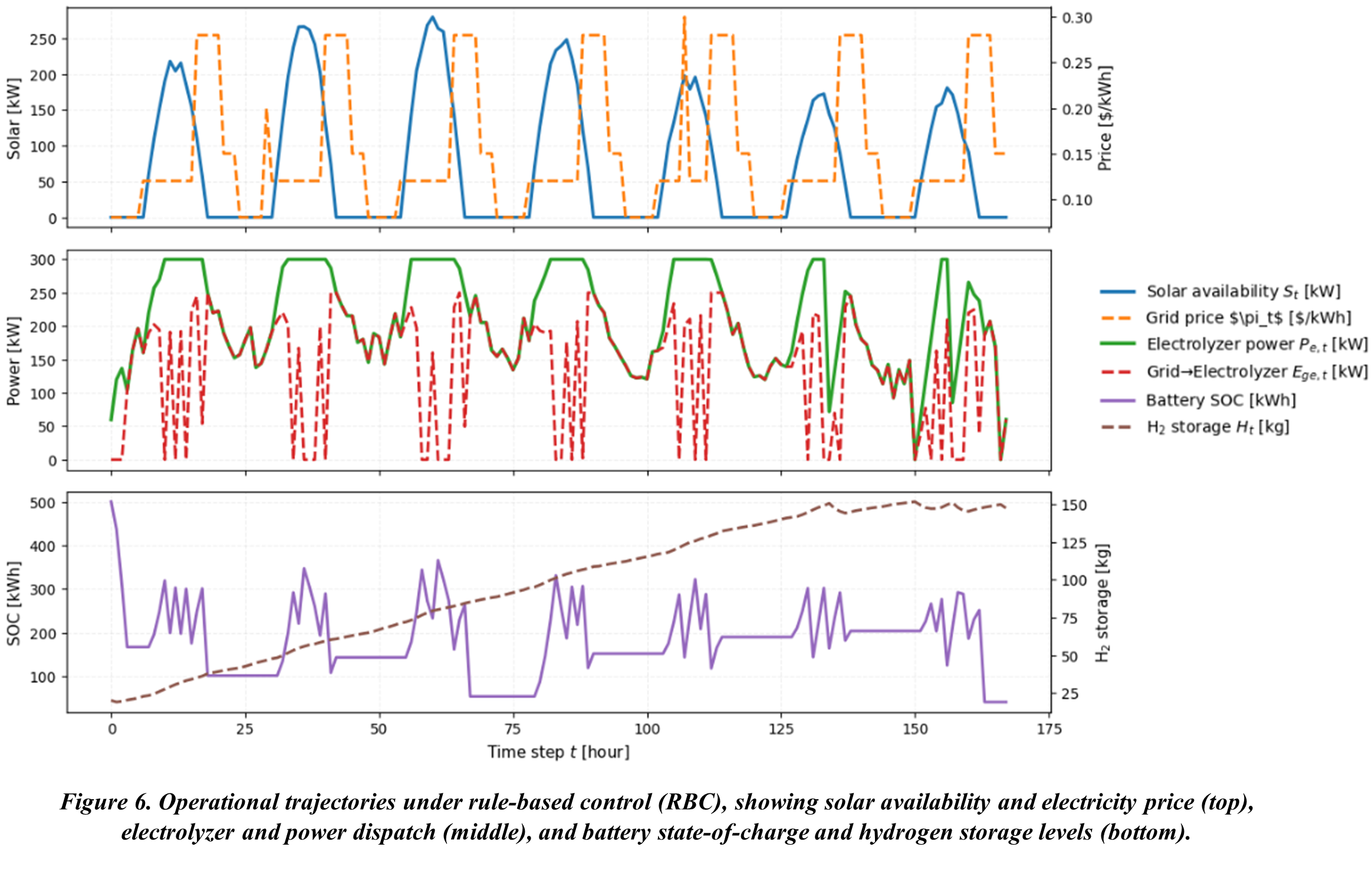}
  \caption{Operational trajectories under rule-based control (RBC): solar
    availability and electricity price (top); electrolyzer and power dispatch
    (middle); battery state of charge and hydrogen storage level (bottom).}
  \label{fig:rbc_trajectories}
\end{figure}

The RBC provides a reactive benchmark without forecasting or optimization.
As shown in Figure~\ref{fig:rbc_trajectories}, solar energy is allocated
instantaneously to the electrolyzer, excess generation charges the battery,
and any remaining shortfall is met by grid imports regardless of price. No
anticipatory charging or production shifting is observed; the battery acts
primarily as a short-term smoothing device, and hydrogen storage is used
passively rather than strategically.

\begin{figure}[htbp]
  \centering
  \includegraphics[width=0.92\textwidth]{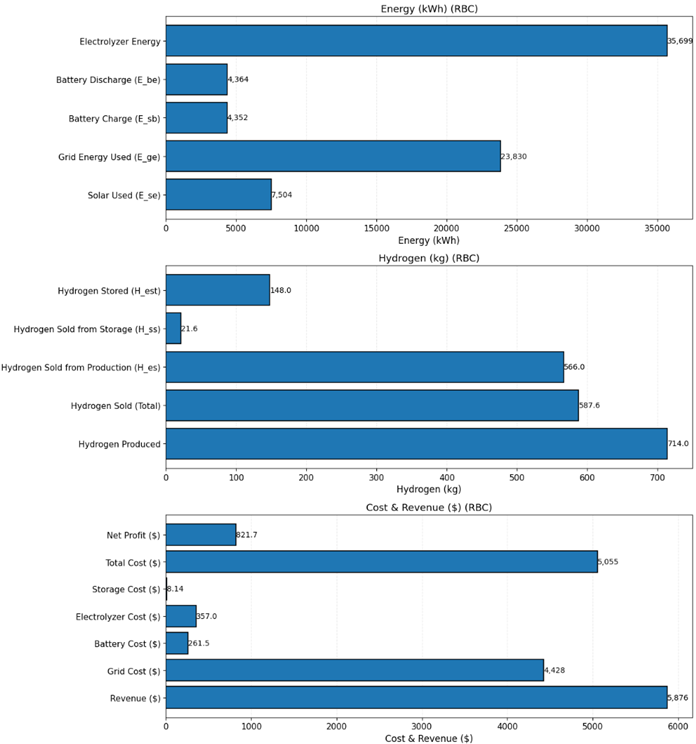}
  \caption{Aggregated performance metrics for rule-based control (RBC): energy
    flows, hydrogen production and utilization, and economic cost--revenue
    breakdown.}
  \label{fig:rbc_metrics}
\end{figure}

The aggregated performance metrics in Figure~\ref{fig:rbc_metrics} reveal three
defining characteristics: high grid dependence; moderate hydrogen utilization
($\sim$80\%), suggesting imperfect alignment between production and demand; and
grid-dominated operating costs that substantially constrain net profitability.
Although RBC maintains stable and feasible operation and achieves positive profit,
it fails to exploit price variability, perform storage arbitrage, or anticipate
future demand shifts. RBC thus establishes a practical lower-bound benchmark
illustrating the structural limitations of reactive dispatch under uncertainty.

\subsection{Model Predictive Control}
\label{subsec:results_mpc}

\begin{figure}[htbp]
  \centering
  \includegraphics[width=0.92\textwidth]{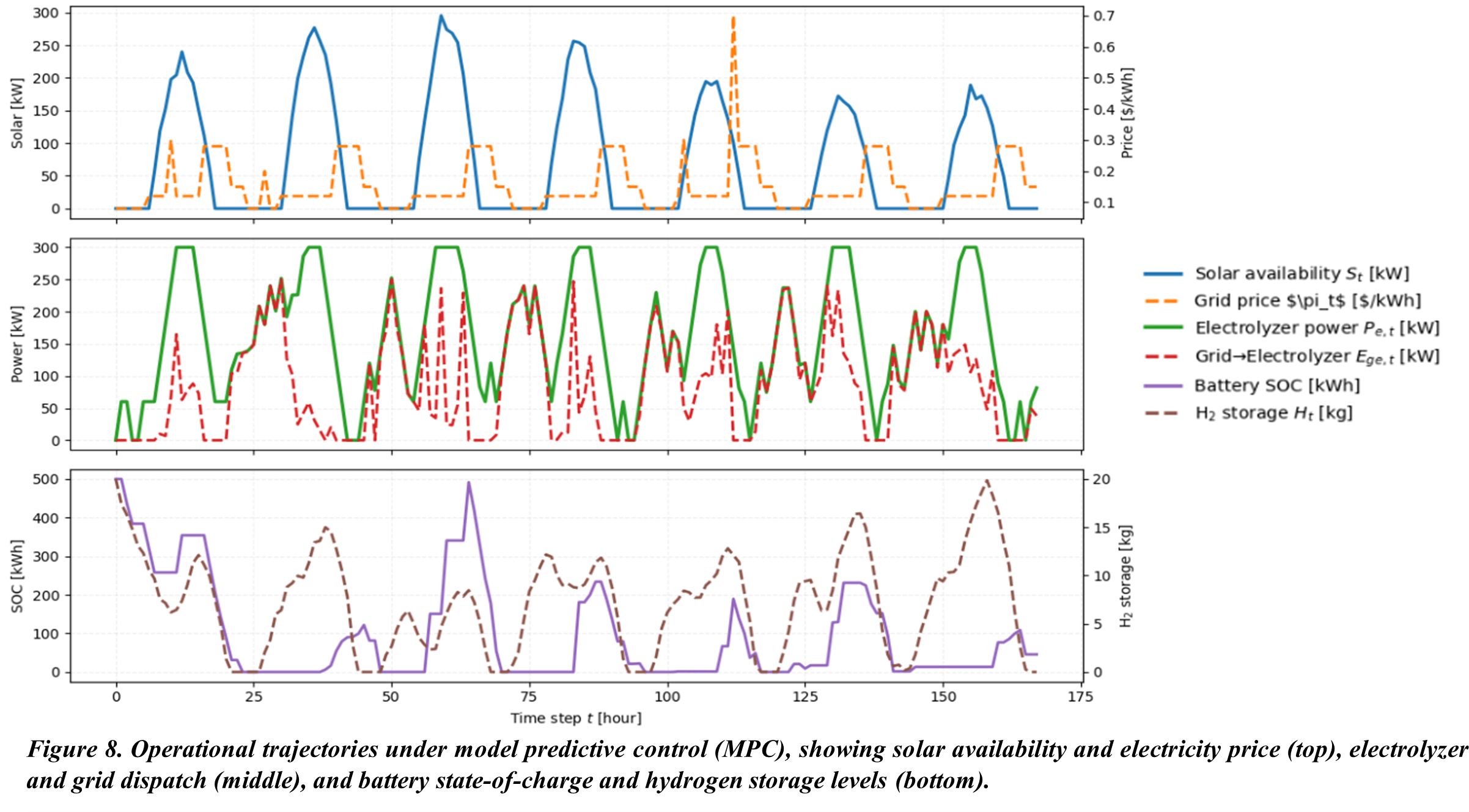}
  \caption{Operational trajectories under model predictive control (MPC): solar
    availability and electricity price (top); electrolyzer and grid dispatch
    (middle); battery state of charge and hydrogen storage level (bottom).}
  \label{fig:mpc_trajectories}
\end{figure}

The operational trajectories in Figure~\ref{fig:mpc_trajectories} show
structured and economically coordinated dispatch behavior. Electrolyzer power
remains consistently high, frequently operating near rated capacity. Grid imports
vary in alignment with price signals, indicating economically rationalized usage
rather than unconditional supplementation. Battery operation is deliberate and
episodic: charging occurs during favorable periods, while discharge is deployed
strategically to support electrolyzer operation during less favorable renewable
conditions. Hydrogen storage exhibits controlled oscillations without persistent
accumulation or depletion, reflecting coordinated temporal management.

\begin{figure}[htbp]
  \centering
  \includegraphics[width=0.92\textwidth]{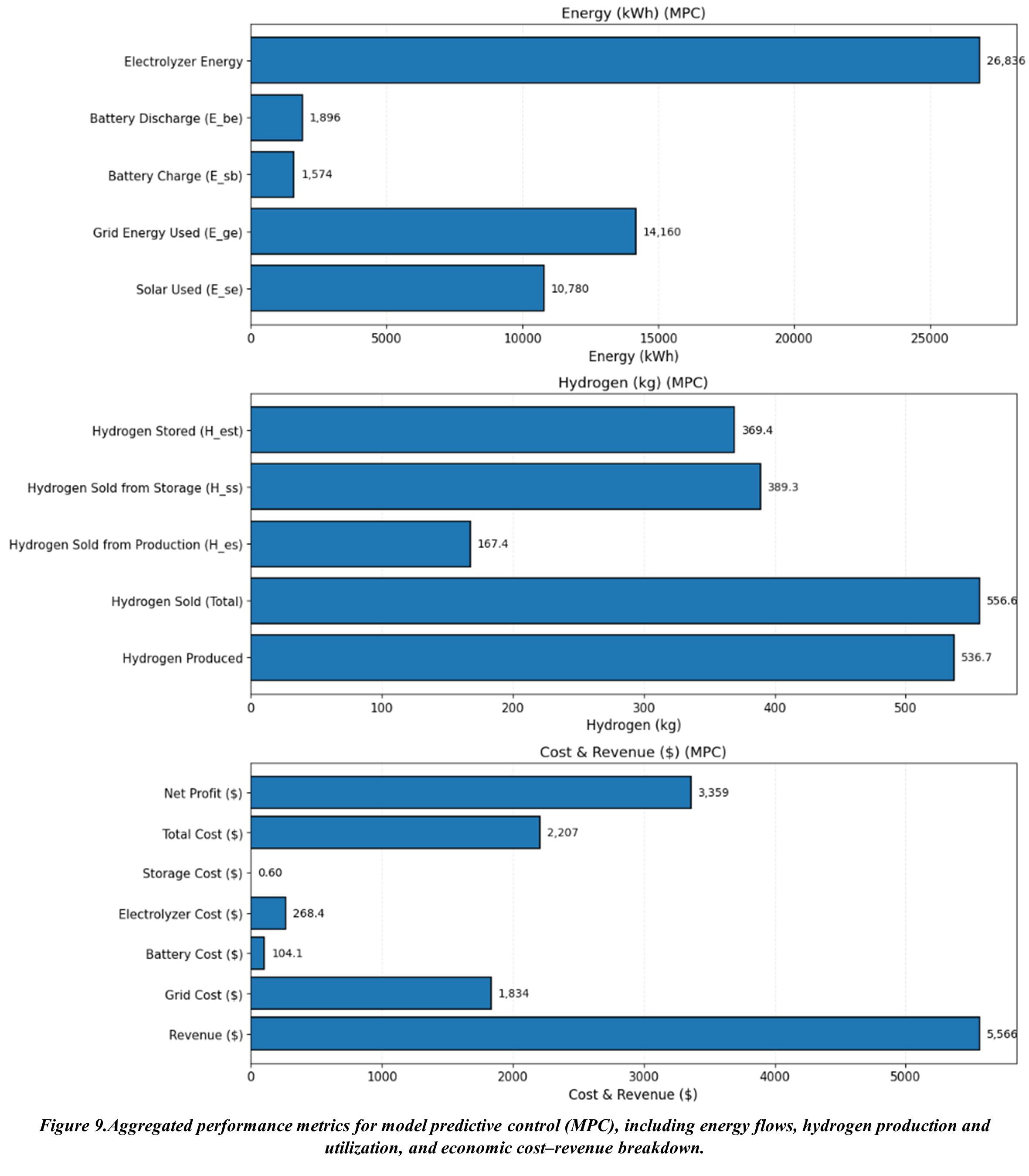}
  \caption{Aggregated performance metrics for model predictive control (MPC):
    energy flows, hydrogen production and utilization, and economic cost--revenue
    breakdown.}
  \label{fig:mpc_metrics}
\end{figure}

The aggregated metrics in Figure~\ref{fig:mpc_metrics} reveal sustained high
electrolyzer utilization, reduced cost intensity relative to throughput, and
strong hydrogen production--sales alignment resulting in near-complete
utilization. Revenue remains dominant relative to total operating cost, and
storage-related costs are negligible. MPC demonstrates stable, economically
efficient operation under realistic constraints by leveraging temporal
coordination across renewable availability, price variability, and demand
realization.

\subsection{Reinforcement Learning without Forecasts (RL-NF)}
\label{subsec:results_rlnf}

\begin{figure}[htbp]
  \centering
  \includegraphics[width=0.92\textwidth]{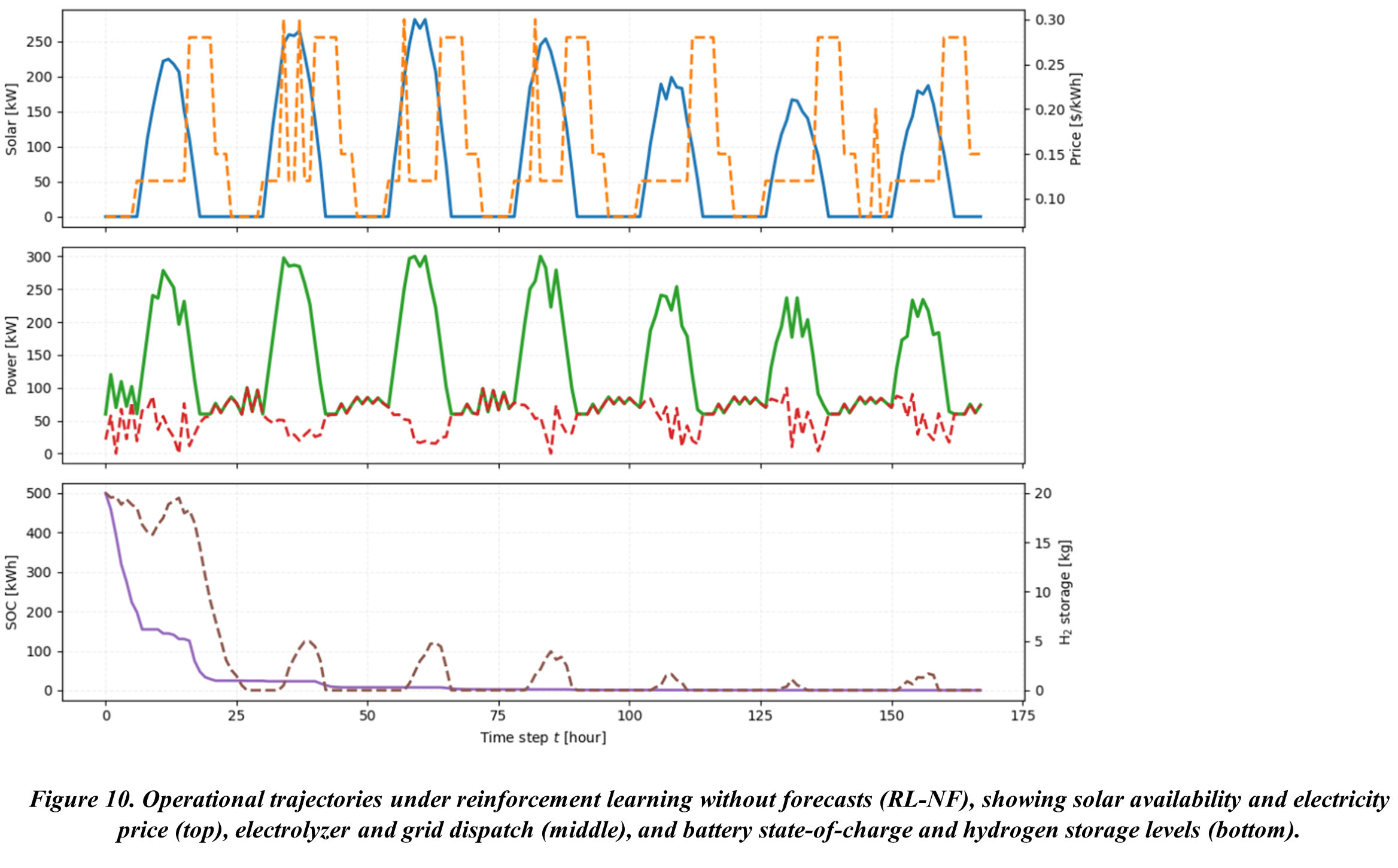}
  \caption{Operational trajectories under reinforcement learning without
    forecasts (RL-NF): solar availability and electricity price (top);
    electrolyzer and grid dispatch (middle); battery state of charge and
    hydrogen storage level (bottom).}
  \label{fig:rlnf_trajectories}
\end{figure}

Under partial observability, RL-NF shows stable but reactive dispatch
(Figure~\ref{fig:rlnf_trajectories}). Electrolyzer power closely tracks
real-time solar availability. Battery participation is minimal, temporal
arbitrage is largely absent, and operational flexibility is achieved primarily
through direct solar--grid balancing. Hydrogen storage is drawn down early and
thereafter remains near low levels.

\begin{figure}[htbp]
  \centering
  \includegraphics[width=0.92\textwidth]{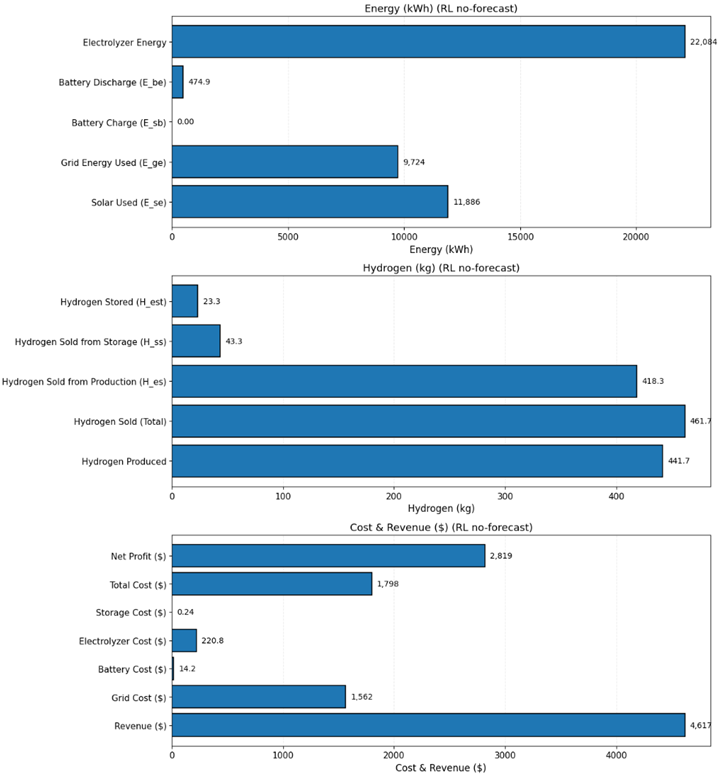}
  \caption{Aggregated performance metrics for reinforcement learning without
    forecasts (RL-NF): energy flows, hydrogen production and utilization, and
    economic cost--revenue breakdown.}
  \label{fig:rlnf_metrics}
\end{figure}

The aggregated results in Figure~\ref{fig:rlnf_metrics} highlight strong direct
solar utilization with moderate grid supplementation, limited battery engagement,
and high production--sales alignment resulting in near-complete hydrogen
utilization. While net profit remains strong, total production scale is lower
than under MPC, reflecting the absence of forecast-driven ramping during
favorable conditions. RL-NF demonstrates that learning-based control can achieve
economically competitive and constraint-respecting operation without explicit
forecasts, though the dispatch pattern remains reactive rather than
intertemporally structured.

\subsection{Reinforcement Learning with Forecasts (RL-F)}
\label{subsec:results_rlf}

\begin{figure}[htbp]
  \centering
  \includegraphics[width=0.92\textwidth]{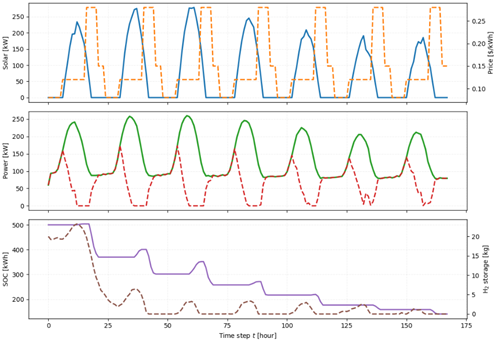}
  \caption{Operational trajectories under reinforcement learning with forecasts
    (RL-F): solar availability and electricity price (top); electrolyzer and grid
    dispatch (middle); battery state of charge and hydrogen storage level
    (bottom).}
  \label{fig:rlf_trajectories}
\end{figure}

Compared to RL-NF, RL-F shows more structured and forward-looking dispatch
(Figure~\ref{fig:rlf_trajectories}). Electrolyzer power does not simply follow
real-time solar availability; operation remains sustained during favorable periods,
indicating the use of short-horizon forecasts for anticipatory scheduling. Grid
imports are adjusted in a coordinated manner. Battery participation remains limited
but purposeful. Hydrogen storage is managed more smoothly, avoiding rapid depletion
and reflecting improved synchronization between production and expected demand.

\begin{figure}[htbp]
  \centering
  \includegraphics[width=0.92\textwidth]{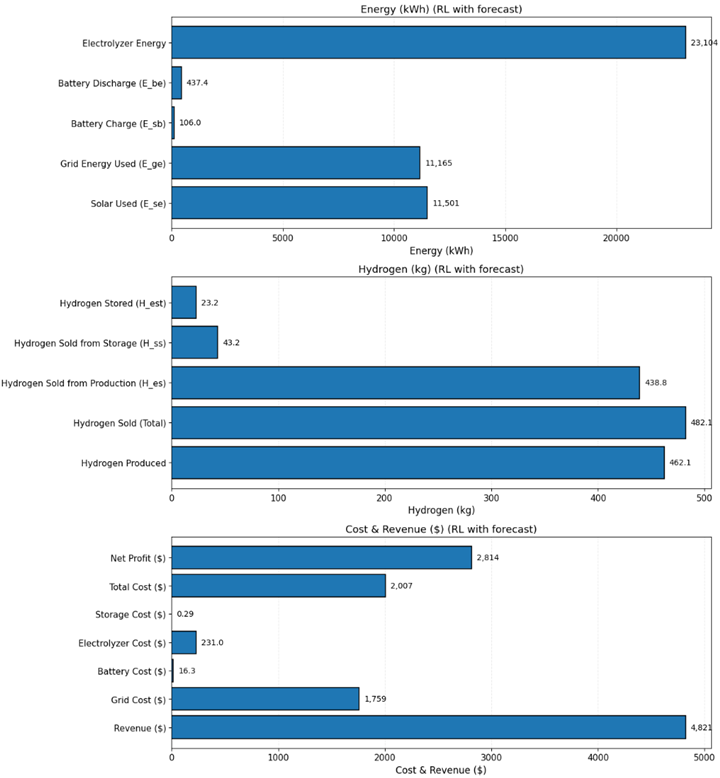}
  \caption{Aggregated performance metrics for reinforcement learning with
    forecasts (RL-F): energy flows, hydrogen production and utilization, and
    economic cost--revenue breakdown.}
  \label{fig:rlf_metrics}
\end{figure}

The aggregated results in Figure~\ref{fig:rlf_metrics} show strong renewable
utilization with balanced grid support, near-complete hydrogen utilization, and
negligible storage and battery costs relative to revenue. Net profit approaches
the MPC benchmark. Compared to RL-NF, total production is similar; the primary
improvement lies in better timing of decisions rather than increased scale.
Incorporating forecasts enhances temporal coordination without requiring explicit
optimization.

\section{Conclusion}
\label{sec:conclusion}

This study presented a unified and physically consistent comparison of four
control strategies for a renewable-powered hydrogen supply chain: RBC, MPC,
RL-NF, and RL-F. All controllers were evaluated under identical disturbances,
constraints, and economic assumptions to ensure a transparent and fair comparison.

Three central findings emerge. First, reactive heuristic control (RBC) maintains
feasibility and positive profitability but exhibits clear structural
inefficiencies, including high grid dependence and reduced hydrogen utilization
($\sim$80\%), highlighting the economic limitations of myopic dispatch under
uncertainty. Second, MPC achieves the strongest economic performance and
near-perfect hydrogen utilization by coordinating electrolyzer operation, storage
dynamics, and grid interaction through explicit forecast-based optimization,
demonstrating the value of structured intertemporal planning when predictive
information is available. Third, RL exhibits competitive and robust behavior,
particularly in the no-forecast setting, achieving high utilization and
substantial profit improvement over heuristic control; however, incorporating
noisy forecasts does not consistently enhance RL performance, indicating that
learning-based exploitation of predictive signals remains condition-dependent.

The analysis relies on several simplifying assumptions, including stylized
disturbance models, a simplified market structure, and a merchant-operation
objective without contractual delivery penalties. RL outcomes are also specific
to the chosen modeling, training, and forecast configurations. Broader validation
under more realistic system representations and multiple uncertainty realizations
is necessary to establish fully general conclusions.

\subsection*{Future Work}

\paragraph{Hybrid MPC--RL.}
Combining MPC's structured, forecast-based optimization with RL's adaptability is
a promising direction. Hybrid designs — such as RL-tuned MPC, MPC-guided RL
training, or supervisory switching frameworks — could preserve constraint
satisfaction and economic efficiency while improving robustness to model mismatch
and changing operating conditions.

\paragraph{Robust MPC.}
Future work should incorporate explicit uncertainty handling via scenario-based,
chance-constrained, or distributionally robust MPC formulations to reduce
sensitivity to forecast errors and improve reliability under renewable and demand
variability.

\paragraph{Multi-Agent Reinforcement Learning (MARL).}
As hydrogen systems scale and decentralize, MARL can enable coordinated control
of distributed assets (electrolyzers, storage units, grid interfaces), supporting
scalability and modular deployment, though stability, coordination, and safety
constraints remain open challenges.

\bibliographystyle{plainnat}
\bibliography{references}

@article{Behzadi2023,
  author  = {Behzadi, A. and Alirahmi, S. M. and Yu, H. and Sadrizadeh, S.},
  title   = {An efficient renewable hybridization based on hydrogen storage for
             peak demand reduction: A rule-based energy control and optimization
             using machine learning techniques},
  journal = {Journal of Energy Storage},
  year    = {2023},
  volume  = {57},
  pages   = {106168},
  doi     = {10.1016/j.est.2022.106168}
}

@article{Brka2016,
  author  = {Brka, A. and Al-Abdeli, Y. M. and Kothapalli, G.},
  title   = {Predictive power management strategies for stand-alone hydrogen
             systems: Operational impact},
  journal = {International Journal of Hydrogen Energy},
  year    = {2016},
  volume  = {41},
  number  = {16},
  pages   = {6685--6698},
  doi     = {10.1016/j.ijhydene.2016.01.031}
}

@article{Glavic2017,
  author  = {Glavic, M. and Fonteneau, R. and Ernst, D.},
  title   = {Reinforcement learning for electric power system decision and
             control: Past considerations and perspectives},
  journal = {IFAC-PapersOnLine},
  year    = {2017},
  volume  = {50},
  number  = {1},
  pages   = {6918--6927},
  doi     = {10.1016/j.ifacol.2017.08.1217}
}

@article{Huang2022,
  author  = {Huang, C. and Zong, Y. and You, S. and Tr{\ae}holt, C.},
  title   = {Economic model predictive control for multi-energy system
             considering hydrogen-thermal-electric dynamics and waste heat
             recovery of {MW}-level alkaline electrolyzer},
  journal = {Energy Conversion and Management},
  year    = {2022},
  volume  = {265},
  pages   = {115697},
  doi     = {10.1016/j.enconman.2022.115697}
}

@book{IEA2019,
  author    = {{International Energy Agency (IEA)}},
  title     = {The Future of Hydrogen: Seizing Today's Opportunities},
  year      = {2019},
  publisher = {International Energy Agency},
  address   = {Paris},
  url       = {https://www.iea.org/reports/the-future-of-hydrogen}
}

@article{Liang2024,
  author  = {Liang, T. and Chai, L. and Cao, X. and Tan, J. and Jing, Y.
             and Lv, L.},
  title   = {Real-time optimization of large-scale hydrogen production systems
             using off-grid renewable energy: Scheduling strategy based on deep
             reinforcement learning},
  journal = {Renewable Energy},
  year    = {2024},
  volume  = {224},
  pages   = {120177},
  doi     = {10.1016/j.renene.2024.120177}
}

@article{Oh2020,
  author  = {Oh, E. and Wang, H.},
  title   = {Reinforcement-learning-based energy storage system operation
             strategies to manage wind power forecast uncertainty},
  journal = {IEEE Access},
  year    = {2020},
  volume  = {8},
  pages   = {20965--20976},
  doi     = {10.1109/ACCESS.2020.2968841}
}

@article{Oldewurtel2012,
  author  = {Oldewurtel, F. and Parisio, A. and Jones, C. N. and Gyalistras, D.
             and Gwerder, M. and Stauch, V. and Morari, M.},
  title   = {Use of model predictive control and weather forecasts for energy
             efficient building climate control},
  journal = {Energy and Buildings},
  year    = {2012},
  volume  = {45},
  pages   = {15--27},
  doi     = {10.1016/j.enbuild.2011.09.022}
}

@article{Parisio2014,
  author  = {Parisio, A. and Rikos, E. and Glielmo, L.},
  title   = {A model predictive control approach to microgrid operation
             optimization},
  journal = {IEEE Transactions on Control Systems Technology},
  year    = {2014},
  volume  = {22},
  number  = {5},
  pages   = {1813--1827},
  doi     = {10.1109/TCST.2013.2295737}
}

@book{Rawlings2020,
  author    = {Rawlings, J. B. and Mayne, D. Q. and Diehl, M. M.},
  title     = {Model Predictive Control: Theory, Computation, and Design},
  year      = {2020},
  edition   = {2nd},
  publisher = {Nob Hill Publishing},
  address   = {Madison, WI}
}

@article{Ruelens2016,
  author  = {Ruelens, F. and Claessens, B. J. and Vandael, S. and
             De~Schutter, B. and Bab\v{u}ska, R. and Belmans, R.},
  title   = {Residential demand response of thermostatically controlled loads
             using batch reinforcement learning},
  journal = {IEEE Transactions on Smart Grid},
  year    = {2016},
  volume  = {8},
  number  = {5},
  pages   = {2149--2159},
  doi     = {10.1109/TSG.2016.2517211}
}

@article{Staffell2019,
  author  = {Staffell, I. and Scamman, D. and Abad, A. V. and Balcombe, P.
             and Dodds, P. E. and Ekins, P. and Schmidt, P. and Bauen, A.
             and Ward, K. R.},
  title   = {The role of hydrogen and fuel cells in the global energy system},
  journal = {Energy \& Environmental Science},
  year    = {2019},
  volume  = {12},
  number  = {2},
  pages   = {463--491},
  doi     = {10.1039/C8EE01157E}
}

@article{VazquezCanteli2019,
  author  = {V{\'a}zquez-Canteli, J. R. and Nagy, Z.},
  title   = {Reinforcement learning for demand response: A review of algorithms
             and modeling techniques},
  journal = {Applied Energy},
  year    = {2019},
  volume  = {235},
  pages   = {1072--1089},
  doi     = {10.1016/j.apenergy.2018.11.002}
}

@article{Velarde2017,
  author  = {Velarde, P. and Valverde, L. and Maestre, J. M. and
             Ocampo-Mart{\'i}nez, C. and Bordons, C.},
  title   = {On the comparison of stochastic model predictive control
             strategies applied to a hydrogen-based microgrid},
  journal = {Journal of Power Sources},
  year    = {2017},
  volume  = {343},
  pages   = {161--173},
  doi     = {10.1016/j.jpowsour.2017.01.015}
}

@article{Yang2021,
  author  = {Yang, T. and Zhao, L. and Li, W. and Zomaya, A. Y.},
  title   = {Dynamic energy dispatch strategy for integrated energy system
             based on improved deep reinforcement learning},
  journal = {Energy},
  year    = {2021},
  volume  = {235},
  pages   = {121377},
  doi     = {10.1016/j.energy.2021.121377}
}

@article{Ye2024,
  author  = {Ye, J. and Wang, X. and Hua, Q. and Sun, L.},
  title   = {Deep reinforcement learning based energy management of a hybrid
             electricity-heat-hydrogen energy system with demand response},
  journal = {Energy},
  year    = {2024},
  volume  = {305},
  pages   = {131874},
  doi     = {10.1016/j.energy.2024.131874}
}

\end{document}